\documentclass[12pt,leqno]{article}

\usepackage{amsmath,amssymb,amsthm}

\makeatletter

\newtheorem{Theorem}{Theorem}

\theoremstyle{definition}

\theoremstyle{remark}

\def\({{\rm (}}
\def\){{\rm )}}

\let\Mathrm\operator@font

\def\standop#1{\mathop{\Mathrm #1}\nolimits}
\def\difstop#1#2{\expandafter\def\csname #1\endcsname{\standop{#2}}}
\def\defstop#1{\difstop{#1}{#1}}

\defstop{Hom}
\defstop{Gal}
\defstop{Spec}

\defstop{Ker}
\defstop{Min}
\difstop{Image}{Im}
\difstop{tdeg}{trans.deg}
\def\red{_{\Mathrm{red}}}
\defstop{length}
\defstop{supp}
\defstop{Proj}
\defstop{Flat}
\difstop{height}{ht}

\def\specialarrow#1{\setbox\z@=\hbox{$\m@th
 \mathop{\vphantom{\rightarrow}}\limits^{\hspace{.5ex}{#1}\hspace
{.8ex}}$}\mathrel{\ifdim\wd\z@<1.2em\dimen\tw@
1.2em\else\dimen\tw@\wd\z@\fi\copy\z@\kern-\wd\z@\hbox to\dimen\tw@
{\rightarrowfill}}}

\def\sdarrow#1{\downarrow\hbox to 0pt{\scriptsize$#1$\hss}}
\def\suarrow#1{\uparrow\hbox to 0pt{\scriptsize$#1$\hss}}

\def\section{\@startsection{section}{1}{\z@ }%
{-3.5ex plus -1ex minus -.2ex}{2.3ex plus .2ex}{\bf }}

\long\def\refname{\par\kern -3ex
\begin{center}\rm R\sc{eferences}\end{center}\par\kern 
-2ex}

\def\@seccntformat#1{\csname the#1\endcsname.\quad}

\def\@@@sect#1#2#3#4#5#6[#7]#8{%
   \ifnum #2>\c@secnumdepth 
      \def \@svsec {}\else \refstepcounter {#1}%
      \def\@svsec{}
   \fi 
   \@tempskipa #5\relax 
   \ifdim \@tempskipa >\z@ 
     \begingroup #6\relax \@hangfrom {\hskip #3\relax 
     \@svsec}{\interlinepenalty \@M #8\par }\endgroup 
     \csname #1mark\endcsname {#7}
   \else 
   \def \@svsechd {#6\hskip #3\@svsec #8\csname #1mark\endcsname {#7}}
   \fi \@xsect {#5}}

\def\@@@startsection#1#2#3#4#5#6{%
 \if@noskipsec \leavevmode \fi \par \@tempskipa #4\relax \@afterindenttrue 
 \ifdim \@tempskipa <\z@ \@tempskipa -\@tempskipa \@afterindentfalse 
 \fi \if@nobreak \everypar {}\else \addpenalty {\@secpenalty }\addvspace 
  {\@tempskipa }\fi \@ifstar {\@ssect {#3}{#4}{#5}{#6}}{\@dblarg 
  {\@@@sect {#1}{#2}{#3}{#4}{#5}{#6}}}}

\def\theparagraph{\thesection.\arabic{paragraph}}
\def\aparagraph{\@@@startsection{paragraph}{2}{\z@ }%
              {1.75ex plus .2ex minus .15ex}{-1em}{\bf(\theparagraph) } }
\def\paragraph{\@@@startsection{paragraph}{2}{\z@ }%
              {1.75ex plus .2ex minus .15ex}{-1em}{}{\bf(\theparagraph)} }

\let\c@Theorem\c@paragraph

\title{A pure subalgebra of a finitely generated algebra 
is finitely generated}

\author{M{\sc itsuyasu} H{\sc ashimoto}}
\date{\normalsize
Graduate School of Mathematics, Nagoya University\\
Chikusa-ku,  Nagoya 464--8602 JAPAN\\
\small \tt hasimoto@math.nagoya-u.ac.jp}

\begin{document}

\maketitle

\begin{abstract}
We prove the following.
Let $R$ be a Noetherian commutative ring, $B$ a finitely generated 
$R$-algebra, and $A$ a pure $R$-subalgebra of $B$.
Then $A$ is finitely generated over $R$.
\end{abstract}

In this paper, all rings are commutative.
Let $A$ be a ring and $B$ an $A$-algebra.
We say that $A\rightarrow B$ is pure, or $A$ is a pure subring of $B$, 
if for any $A$-module $M$, the map $M=M\otimes_A A \rightarrow M\otimes_A B$ 
is injective.
Considering the case $M=A/I$, where $I$ is an ideal of $A$, 
we immediately have that $IB\cap A=I$.

It has been shown that if $B$ has a good property and $A$ is a pure
subring of $B$, then $A$ has a good property.
If $B$ is a regular Noetherian ring containing a field, then $A$ is 
Cohen-Macaulay \cite{HR}, \cite{HH}.
If $k$ is a field of characteristic zero, $A$ and $B$ are 
essentially of finite type over $k$, and $B$ has at most 
rational singularities, then $A$ has at most rational singularities 
\cite{Boutot}.

In this paper, we prove the following

\begin{Theorem}\label{main.thm}
Let $R$ be a Noetherian ring, 
$B$ a finitely generated $R$-algebra, and $A$ a pure $R$-subalgebra of $B$.
Then $A$ is finitely generated over $R$.
\end{Theorem}

The case that $B$ is $A$-flat is proved in \cite{Hashimoto}.
This theorem is on the same line as the finite generation 
results in \cite{Hashimoto}.

To prove the theorem, we need the following, which is a special case of 
a theorem of Raynaud-Gruson \cite{Raynaud}, \cite{RG}.

\begin{Theorem}\label{flattening.thm}
Let $A\rightarrow B$ be a homomorphism of Notherian rings, and 
$\varphi\colon X\rightarrow Y$ the associated morphism of affine schemes.
Let $U\subset Y$ be an open subset, and assume that $\varphi\colon
\varphi^{-1}(U)\rightarrow U$ is flat.
Then there exists some ideal $I$ of $A$ such that $V(I)\cap U=\emptyset$, 
and that the morphism $\Phi\colon \Proj R_B(BI)\rightarrow \Proj R_A(I)$, 
determined by the associated morphism of the Rees algebras 
$R_A(I):=A[tI]\rightarrow R_B(BI):=B[tBI]$, is flat.
\end{Theorem}

The morphism $\Phi$ in the theorem is called a flattening of $\varphi$.

\proof[Proof of Theorem~\ref{main.thm}]
Note that for any $A$-algebra $A'$, the homomorphism
$A'\rightarrow B\otimes_A A'$ is pure.

Since $B$ is finitely generated over $R$, it is Noetherian.
Since $A$ is a pure subring of $B$, $A$ is also Noetherian.
So if $A\red$ is finitely generated, then so is $A$.
Replacing $A$ by $A\red$ and $B$ by $B\otimes_A A\red$, we may 
assume that $A$ is reduced.

Since $A\rightarrow \prod_{P\in\Min(A)}A/P$ is finite and injective,
it suffices to prove that each $A/P$ is finitely generated for 
$P\in\Min(A)$, where $\Min(A)$ denotes the set of minimal primes of $A$.
By the base change, we may assume that $A$ is a domain.

There exists some minimal prime $P$ of $B$ such that $P\cap A=0$.
Assume the contrary.
Then take $a_P\in P\cap A\setminus\{0\}$ for each $P\in\Min(B)$.
Then $\prod_P a_P$ must be nilpotent, which contradicts to our 
assumption that $A$ is a domain.

So by \cite[(2.11) and (2.20)]{Onoda}, $A$ is a finitely generated 
$R$-algebra if and only if $A_{\mathfrak p}$ is a finitely generated 
$R_{\mathfrak p}$-algebra for each $\mathfrak p\in\Spec R$.
So we may assume that $R$ is a local ring.

By the descent argument \cite[(2.7.1)]{EGA-IV}, $\hat R\otimes_R A$ is 
a finitely generated $\hat R$-algebra if and only if $A$ is a finitely 
generated $R$-algebra, where $\hat R$ is the completion of $R$.
So we may assume that $R$ is a complete local ring.
We may lose the assumption that $A$ is a domain (even if $A$ is a 
domain, $\hat R\otimes_R A$ may not be a domain).
However, doing the same reduction argument as above if necessary, 
we may still assume that $A$ is a domain.

Let $\varphi\colon X\rightarrow Y$ be a morphism of affine schemes
associated with the map $A\rightarrow B$.
Note that $\varphi$ is a morphism of finite type between Noetherian
schemes.
We denote the flat locus of $\varphi$ by $\Flat(\varphi)$.
Then $\varphi(X\setminus\Flat(\varphi))$ is a constructible set of $Y$ 
not containing the generic point.
So $U=Y\setminus\overline{\varphi(X\setminus\Flat(\varphi))}$ is a dense
open subset of $Y$, and $\varphi\colon \varphi^{-1}(U)\rightarrow U$ is
flat.
By Theorem~\ref{flattening.thm}, there exists some nonzero ideal $I$ of 
$A$ such that $\Phi\colon \Proj R_B(BI)\rightarrow \Proj R_A(I)$ is flat.

If $J$ is a homogeneous ideal of $R_A(I)$, then we have an expression 
$J=\bigoplus_{n\geq 0} J_nt^n$ ($J_n\subset I^n$).
Since $A$ is a pure subalgebra of $B$, we have $J_n B\cap I^n=J_n$ for 
each $n$.
Since $J R_B(BI)=\bigoplus_{n\geq 0}(J_nB)t^n$, we have that 
$J R_B(BI)\cap R_A(I)= J$.
Namely, any homogeneous ideal of $R_A(I)$ is contracted from $R_B(BI)$.

Let $P$ be a homogeneous prime ideal of $R_A(I)$.
Then there exists some minimal prime $Q$ of $PR_B(BI)$ such that 
$Q\cap R_A(AI)=P$.
Assume the contrary.
Then for each minimal prime $Q$ of $PR_B(BI)$, there exists some 
$a_Q\in (Q\cap R_A(AI))\setminus P$.
Then $\prod a_Q\in \sqrt{PR_B(BI)}\cap R_A(AI)\setminus P$.
However, we have
\[
\sqrt{P R_B(BI)}\cap R_A(I)=\sqrt{P R_B(BI)\cap R_A(I)}=\sqrt{P}=P,
\]
and this is a contradiction.
Hence $\Phi\colon \Proj R_B(BI)\rightarrow \Proj R_A(I)$ is 
faithfully flat.

Since $\Proj R_B(BI)$ is of finite type over $R$ and $\Phi$ is faithfully
flat, we have that $\Proj R_A(I)$ is of finite type by 
\cite[Corollary~2.6]{Hashimoto}.
Note that the blow-up $\Proj R_A(I)\rightarrow Y$ is proper surjective.
Since $R$ is excellent, $Y$ is of finite type over $R$ by
\cite[Theorem~4.2]{Hashimoto}.
Namely, $A$ is a finitely generated $R$-algebra.
\qed

\end{document}